\documentclass[preprint,3p]{elsarticle}
\usepackage{amsmath}
\usepackage{amssymb}
\usepackage{caption} 
\usepackage{xcolor}
\usepackage{color}
 \usepackage[english]{babel} 
 \usepackage{float} 
 \usepackage{graphics} 
 \usepackage{graphicx} 
\usepackage[utf8]{inputenc}
\usepackage{mathtools}
\usepackage{capt-of}
\usepackage{calrsfs}
\usepackage{lipsum}  
 \usepackage{epsfig}
\newtheorem{theorem}{Theorem}[section]
\newtheorem{lemma}{Lemma}[section]

\newtheorem{definition}{Definition}[section]

\newtheorem{pro}{Property}[section]
\newdefinition{remark}{Remark}[section]
\newdefinition{fe}{Feature}[section]
\newproof{proof}{Proof}
\newproof{pot}{Proof of Theorem \ref{thm2}}
\newproof{poot}{Proof of Corollary \ref{co1}}
\numberwithin{equation}{section}
\newdefinition{exe}{Example}[section]
\journal{}
\begin{document}
\begin{frontmatter}

\title{Fractional Sturm-Liouville eigenvalue problems, II }
\author{
M. Dehghan\footnote[1]{dehghan@math.carleton.ca} and A. B. Mingarelli\footnote[1]{angelo@math.carleton.ca}}
\address{$^{1}$School of Mathematics and Statistics, Carleton University, Ottawa, Canada}

\begin{abstract}
We continue the study of a non self-adjoint fractional three-term Sturm-Liouville boundary value problem (with a potential term)  formed by the composition of a left Caputo and left-Riemann-Liouville fractional integral under {\it Dirichlet type} boundary conditions. We study the existence and asymptotic behavior of the real eigenvalues and show that for certain values of the fractional differentiation parameter $\alpha$, $0<\alpha<1$, there is a finite set of real eigenvalues and that, for $\alpha$ near $1/2$, there may be none at all. As $\alpha \to 1^-$ we show that their number becomes infinite and that the problem then approaches a standard Dirichlet Sturm-Liouville problem with the composition of the operators becoming the operator of second order differentiation.
\end{abstract}

\end{frontmatter}

\section{Introduction}

This is a continuation of \cite{dm1} where the results therein are extended to three-term Fractional Sturm-Liouville operators (with a potential term)  formed by the composition of a left Caputo and left-Riemann-Liouville fractional integral. Specifically, the boundary value problem is of the form,
\begin{equation}\label{fe30}
-^{c}\mathcal{D}_{0^+}^{\alpha}\circ   \mathcal{D}^{\alpha}_{0^+}y(t)+q(t)y(t)=\lambda y(t),\qquad 1/2 <\alpha<1,\quad 0 \leq t \leq 1,
\end{equation}
with boundary conditions 
\begin{equation}\label{fcon0}
\mathcal{I}_{0^+}^{1-\alpha}y(t)|_{t=0}=c_1,\quad \text{and}\quad \mathcal{I}_{0^+}^{1-\alpha}y(t)|_{t=1}=c_2,
\end{equation}
where $c_1,c_2$ are real constants and the real valued unspecified potential function, $q\in L^\infty[0,1]$. We note that these are not self-adjoint problems and so there may be non-real spectrum, in general. A well-known property of the Riemann-Liouville integral gives that if the solutions are continuous on $[0,1]$ then the boundary conditions \eqref{fcon0} reduce to the usual fixed end boundary conditions, $y(0)=y(1)=0$, as $\alpha \to 1$. 

For the analogue of the Dirichlet problem described above we study the existence and asymptotic behavior of the real eigenvalues and show that for each $\alpha$, $0<\alpha<1$, there is a finite set of real eigenvalues and that, for $\alpha$ near $1/2$, there may be none at all. As $\alpha \to 1^-$ we show that their number becomes infinite and that the problem then approaches a standard Dirichlet Sturm-Liouville problem with the composition of the operators becoming the operator of second order differentiation acting on a suitable function space.

\section{Preliminaries}  
We recall some definitions from Fractional Calculus and refer the reader to our previous paper \cite{dm1} for further details.

\begin{definition}
The left and the right Riemann-Liouville fractional integrals $\mathcal{I}^{\alpha}_{a^+}$ and $\mathcal{I}^{\alpha}_{b^-}$ of order $\alpha\in \mathbb{R}^+$ are defined by
\begin{equation}\label{lfi}
\mathcal{I}^{\alpha}_{a^+}f(t):=\frac{1}{\Gamma(\alpha)}\int_a^t\frac{f(s)}{(t-s)^{1-\alpha}}ds,\quad t \in (a,b],
\end{equation}
and
\begin{equation}\label{rfi}
\mathcal{I}^{\alpha}_{b^-}f(t):=\frac{1}{\Gamma(\alpha)}\int_t^b\frac{f(s)}{(s-t)^{1-\alpha}}ds,\quad t \in [a,b),
\end{equation}
respectively. Here $\Gamma (\alpha)$ denotes  Euler's Gamma function. The following property is easily verified.
 \end{definition}
 \begin{pro}\label{c}
\normalfont 
 For a constant $C$, we have $\mathcal{I}_{a^+}^{\alpha}C=\frac{(t-a)^{\alpha}}{\Gamma(\alpha+1)}\cdot C$.
 \end{pro}
 The proof is by direct calculation.

  \begin{definition}
The left and the right Caputo fractional derivatives $^{c}\mathcal{D}^{\alpha}_{a^+}$ and $^{c}\mathcal{D}^{\alpha}_{b^-}$ are defined by
\begin{equation}\label{lcfd}
^{c}\mathcal{D}^{\alpha}_{a^+}f(t):=\mathcal{I}^{n-\alpha}_{a^+}\circ\mathcal{D}^nf(t)=\frac{1}{\Gamma(n-\alpha)}\int_a^{t}\frac{f^{(n)}(s)}{(t-s)^{\alpha-n+1}}ds, \quad t>a,
\end{equation}
and 
\begin{equation}\label{rcfd}
^{c}\mathcal{D}^{\alpha}_{b^-}f(t):=(-1)^n\mathcal{I}^{n-\alpha}_{b^-}\circ\mathcal{D}^nf(t)=\frac{(-1)^n}{\Gamma(n-\alpha)}\int_t^{b}\frac{f^{(n)}(s)}{(s-t)^{\alpha-n+1}}ds, \quad t<b,
\end{equation}
respectively, where $f$ is sufficiently differentiable and $n-1\leq \alpha < n$.
 \end{definition}
  \begin{definition}
Similarly, the left and the right Riemann-Liouville  fractional derivatives $\mathcal{D}^{\alpha}_{a^+}$ and $\mathcal{D}^{\alpha}_{b^-}$ are defined by
\begin{equation}\label{lrlfd}
\mathcal{D}^{\alpha}_{a^+}f(t):=\mathcal{D}^n\circ\mathcal{I}^{n-\alpha}_{a^+}f(t)=\frac{1}{\Gamma(n-\alpha)}\frac{d^n}{dt^n}\int_a^{t}\frac{f(s)}{(t-s)^{\alpha-n+1}}ds, \quad t>a,
\end{equation}
and 
\begin{equation}\label{rrlfd}
\mathcal{D}^{\alpha}_{b^-}f(t):=(-1)^n\mathcal{D}^n\circ\mathcal{I}^{n-\alpha}_{b^-}f(t)=\frac{(-1)^n}{\Gamma(n-\alpha)}\frac{d^n}{dt^n}\int_t^{b}\frac{f(s)}{(s-t)^{\alpha-n+1}}ds, \quad t<b,
\end{equation}
respectively, where $f$ is sufficiently differentiable and $n-1\leq \alpha < n$.
 \end{definition}
 \begin{pro}\label{Dt}
 For $\Re(\nu)>-1$, $0<\alpha<1$, and $t>0$, we have
 \[
 D^{\alpha}_{0^+}(t^{\nu})=\frac{\Gamma(1+\nu)}{\Gamma(1+\nu-\alpha)}t^{\nu-\alpha}
 \]
 \end{pro}
 \begin{pro}\label{Dtc}
 For $\Re(\nu)>0$, $0<\alpha<1$,  and $t>0$, we have
 \[
 ^{c}D^{\alpha}_{0^+}(t^{\nu})=\frac{\Gamma(1+\nu)}{\Gamma(1+\nu-\alpha)}t^{\nu-\alpha}
 \]
 \end{pro}
 \begin{pro}\label{id}
 If $y(t)\in L^1(a,b)$ and $\mathcal{I}^{1-\alpha}_{a^+}y, \mathcal{I}^{1-\alpha}_{b^-}y\in AC[a,b]$, then
 \begin{equation}\nonumber
 \begin{split}
 \mathcal{I}^{\alpha}_{a^+}\mathcal{D}^{\alpha}_{a^+}y(t)&=y(t)-\frac{(t-a)^{\alpha-1}}{\Gamma(\alpha)}\mathcal{I}^{1-\alpha}_{a^+}y(a),\\
 \mathcal{I}^{\alpha}_{b^-}\mathcal{D}^{\alpha}_{b^-}y(t)&=y(t)-\frac{(b-t)^{\alpha-1}}{\Gamma(\alpha)}\mathcal{I}^{1-\alpha}_{b^-}y(b).
 \end{split}
 \end{equation}
 \end{pro}
 
  \begin{pro}\label{id_2}
 If $y(t)\in AC[a,b]$ and $0<\alpha\leq 1$, then
 \begin{equation}\nonumber
 \begin{split}
 \mathcal{I}^{\alpha}_{a^+} {^{c}\mathcal{D}^{\alpha}_{a^+}}y(t)&=y(t)-y(a),\\
 \mathcal{I}^{\alpha}_{b^-} {^{c}\mathcal{D}^{\alpha}_{b^-}}y(t)&=y(t)-y(b).
 \end{split}
 \end{equation}
 \end{pro}
 \begin{pro}\label{ricap}
 For $0<\alpha<1$ we have
 \[
 \mathcal{D}^{\alpha}_{a^+}f(t)=\frac{f(a)}{\Gamma(1-\alpha)}(t-a)^{-\alpha}+ ^{c}\mathcal{D}^{\alpha}_{a^+}f(t)
 \] 
 \end{pro}

 \subsection{The Mittag-Leffler function}
The function $E_{\delta}(z)$ defined by 
 \begin{equation}\label{mittag1}
 E_{\delta}(z):=\sum_{k=0}^{\infty}\frac{z^{\delta}}{\Gamma(\delta k+1)},\quad (z\in \mathbb{C}, \Re(\delta)>0),
 \end{equation}
 was introduced by Mittag-Leffler \cite{mit}. In particular, when $\delta=1$ and $\delta=2$, we have
 \begin{equation}\label{e1}
 E_1(z)=e^z,\qquad E_2(z)=\cosh(\sqrt{z}).
\end{equation}
 The generalized Mittag-Leffler function $E_{\delta,\theta}(z)$ is defined by
\begin{equation}\label{mittag2}
E_{\delta,\theta}(z) = \sum_{k=0}^{\infty}\, \frac{z^k}{\Gamma (\delta k+\theta)},
\end{equation}
where $z,\theta \in \mathbb{C}$ and ${\rm Re}\, (\delta)>0$. When $\theta=1$, $E_{\delta,\theta}(z)$ coincides with the Mittag-Leffler function (\ref{mittag1}):
 \begin{equation}\label{e2}
 E_{\delta,1}(z)=E_{\delta}(z).
 \end{equation}
 Two other particular cases of (\ref{mittag2}) are as follows:
 \begin{equation}\label{pmi}
 E_{1,2}(z)=\frac{e^z-1}{z},\quad E_{2,2}(z)=\frac{\sinh(\sqrt{z})}{\sqrt{z}}.
 \end{equation}
 \begin{pro}\label{rele}
 For any $\delta$ with $\Re(\delta)>0$ and for any $z\neq 0$ we have
 \[
 E_{\delta,\delta}(z)=\frac{1}{z}E_{\delta,0}(z)
 \]
 \end{pro}
Further properties of this special function may be found in  \cite{erdelyi}.

 \begin{pro}\label{asym}\rm{ (See \cite{kil}, p.43.)}
If $0<\delta<2$ and $\mu\in(\frac{\delta \pi}{2},\min(\pi,\delta \pi))$, then function $E_{\delta,\theta}(z)$ has the following exponential expansion as $|z|\rightarrow\infty$
 \begin{gather}
E_{\delta,\theta}(z)=\left\{
\begin{array}{ll}
\frac{1}{\delta}z^{\frac{1-\theta}{\delta}}\exp(z^{\frac{1}{\delta}})-\sum_{k=1}^N\frac{1}{\Gamma(\theta-\delta k)}\frac{1}{z^k}+O(\frac{1}{z^{N+1}}), \qquad |\arg(z)|\leq \mu, \label{asyx1}\\ \\
-\sum_{k=1}^N\frac{1}{\Gamma(\theta-\delta k)}\frac{1}{z^k}+O(\frac{1}{z^{N+1}}), \qquad \mu\leq|\arg(z)|\leq \pi.
\end{array}\right. 
\end{gather}
\end{pro}

\section{Existence and uniqueness of the solution of SLPs}
First, we proceed formally. Separating terms in (\ref{fe30}), we get
\begin{equation}\nonumber
^{c}\mathcal{D}_{0^+}^{\alpha}\circ   \mathcal{D}^{\alpha}_{0^+}y(t)=(q(t)-\lambda) y(t),\qquad 1/2 <\alpha<1,\quad 0 \leq t \leq 1.
\end{equation}
Taking the  left  Riemann-Liouville fractional integrals $\mathcal{I}^{\alpha}_{a^+}$ on both sides of the above equation and using Property \ref{id_2},  we have
\begin{equation}\nonumber
\mathcal{D}_{0^+}^{\alpha}y(t)-
\mathcal{D}_{0^+}^{\alpha}y(t)|_{t=0}=\mathcal{I}_{0^+}^{\alpha}((q(t)-\lambda)y(t)).
\end{equation}
Taking the  left  Riemann-Liouville fractional integrals $\mathcal{I}^{\alpha}_{a^+}$ from both sides of the above equation once again and using Property \ref{id},  we get
\begin{equation}\nonumber
y(t)-\frac{t^{\alpha-1}}{\Gamma(\alpha)}\mathcal{I}^{1-\alpha}_{0^+}y(t)|_{t=0}-\mathcal{I}^{\alpha}_{0^+}(\mathcal{D}_{0^+}^{\alpha}y(t)|_{t=0})=\mathcal{I}^{\alpha}_{0^+}(\mathcal{I}_{0^+}^{\alpha}((q(t)-\lambda)y(t)))
\end{equation}
Using Property \ref{c}, we can write
\begin{equation}\nonumber
y(t) = c_1\frac{t^{\alpha-1}}{\Gamma(\alpha)}+c_2\frac{t^{\alpha}}{\Gamma(\alpha+1)}+\mathcal{I}^{\alpha}_{0^+}(\mathcal{I}_{0^+}^{\alpha}((q(t)-\lambda)y(t)))
\end{equation}
in which
\begin{equation}\nonumber
c_1=\mathcal{I}^{1-\alpha}_{0^+}y(t)|_{t=0},\qquad c_2=\mathcal{D}^{\alpha}_{0^+}y(t)|_{t=0}.
\end{equation}
We get through the double fractional integral in the above equation as follows
\begin{equation}\nonumber
y(t) = c_1\frac{t^{\alpha-1}}{\Gamma(\alpha)}+c_2\frac{t^{\alpha}}{\Gamma(\alpha+1)}+\frac{1}{\Gamma^2(\alpha)}\int_0^t(t-s)^{\alpha-1}\left(\int_0^s\frac{(q(r)-\lambda)y(r)}{(s-r)^{1-\alpha}}dr\right)ds.
\end{equation}
By changing the order of integrals in the above equation we get
\begin{equation}\nonumber
y(t) = c_1\frac{t^{\alpha-1}}{\Gamma(\alpha)}+c_2\frac{t^{\alpha}}{\Gamma(\alpha+1)}+\frac{1}{\Gamma^2(\alpha)}\int_0^t(q(r)-\lambda)y(r)\left(\int_r^t(t-s)^{\alpha-1}(s-r)^{\alpha-1}ds\right)dr
\end{equation}
Solving the inner integral gives us
\begin{equation}\label{IE0}
y(t,\lambda) = c_1\frac{t^{\alpha-1}}{\Gamma(\alpha)}+c_2\frac{t^{\alpha}}{\Gamma(\alpha+1)}+\frac{1}{\Gamma(2\alpha)}\int_0^t(q(s)-\lambda)y(s,\lambda)(t-s)^{2\alpha-1}ds.
\end{equation}
We will now show that \eqref{IE0} has a solution that exists in a neighbourhood of $t=0$ and is unique there. Working backwards will then provide us with a unique solution to \eqref{fe30}-\eqref{fcon0}. Although this result already appears in \cite{msjg}, we give a shorter proof part of which will be required later. \\

To this end, let $t > 0$. Define
\begin{equation}\label{iter}
y_n(t,\lambda) = y_0(t,\lambda)+\frac{1}{\Gamma(2\alpha)}\int_0^t(t-s)^{2\alpha-1}(q(s)-\lambda)y_{n-1}(s,\lambda)ds,
\end{equation}
where
\begin{equation}\label{y_0}
y_0(t,\lambda) = c_1\frac{t^{\alpha-1}}{\Gamma(\alpha)}+c_2\frac{t^{\alpha}}{\Gamma(\alpha+1)}.
\end{equation}
Let $\lambda \in \mathbf{C}$, $|\lambda|< \Lambda$, where $\Lambda>0$ is arbitrary but fixed. Then,
\begin{equation}\label{n_1}
\begin{split}
|y_1(t,\lambda)-y_0(t,\lambda)|&\leq \frac{1}{\Gamma(2\alpha)}\int_0^t(t-s)^{2\alpha-1}|q(s)-\lambda||y_0(s,\lambda)|ds\\
&\leq \frac{||q||_{\infty}+\Lambda}{\Gamma(2\alpha)}\int_0^t(t-s)^{2\alpha-1}|y_0(s,\lambda)|ds,
\end{split}
\end{equation}
in which $||q||_{\infty} =\sup_{t\in[0,1]}|q(t)|$.
Substituting (\ref{y_0}) in (\ref{n_1}) and using the fact that,
\begin{equation}\nonumber
\int_a^t(t-s)^{\alpha-1}(s-a)^{\beta-1}ds=\frac{(t-a)^{\alpha+\beta-1}\Gamma(\alpha)\Gamma(\beta)}{\Gamma(\alpha+\beta)},
\end{equation}
we have
\begin{equation}\nonumber
|y_1(t,\lambda)-y_0(t,\lambda)|\leq(||q||_{\infty}+\Lambda)\left (\frac{c_1}{\Gamma(3\alpha)}t^{3\alpha-1}+c_2\frac{t^{3\alpha}}{\Gamma(3\alpha+1)}\right ),
\end{equation}
Now, for $n=2$ in (\ref{iter}) we get
\begin{equation}\label{n_2}
\begin{split}
|y_2(t,\lambda)-y_1(t,\lambda)|&\leq \frac{1}{\Gamma(2\alpha)}\int_0^t(t-s)^{2\alpha-1}|q(s)-\lambda||y_1(s,\lambda)-y_0(s,\lambda)|ds\\
&\leq \frac{1}{\Gamma(2\alpha)}\int_0^t(t-s)^{2\alpha-1}|q(s)-\lambda|\left((||q||_{\infty}+\Lambda)(\frac{c_1}{\Gamma(3\alpha)}s^{3\alpha-1}+c_2\frac{s^{3\alpha}}{\Gamma(3\alpha+1)})\right)ds\\
&\leq (||q||_{\infty}+\Lambda)^2\left(\frac{c_1}{\Gamma(5\alpha)}t^{5\alpha-1}+\frac{c_2}{\Gamma(5\alpha+1)}t^{5\alpha}\right).
\end{split}
\end{equation}
Continuing in this way we get that the series
\begin{equation}\label{y_succ}
 y_0(t,\lambda) + \sum_{n=1}^{\infty}(y_n(t,\lambda)-y_{n-1}(t,\lambda))
\end{equation}
where
\begin{equation}\label{y_sum}
\begin{split}
\sum_{n=1}^{\infty}|y_n(t,\lambda)-y_{n-1}(t,\lambda)|&\leq c_1t^{-1}\sum_{n=1}^{\infty}\frac{(||q||_{\infty}+\Lambda)^n}{\Gamma(2n\alpha+\alpha)}t^{2n\alpha+\alpha}+c_2\sum_{n=1}^{\infty}
\frac{(||q||_{\infty}+\Lambda)^n}{\Gamma(2n\alpha+\alpha+1)}t^{2n\alpha+\alpha}.
\end{split}
\end{equation}
converges uniformly on compact subsets of $(0,1]$. Denote the sum of the infinite series in \eqref{y_succ} by $y(t, \lambda)$. So, by virtue of (\ref{y_0}) and (\ref{y_sum}),  (\ref{y_succ}) gives us,
\begin{equation}\nonumber
\begin{split}
|y(t,\lambda)| &\leq |y_0(t,\lambda)| + \sum_{n=1}^{\infty}|y_n(t,\lambda)-y_{n-1}(t,\lambda)|\\
& \leq \frac{c_1}{\Gamma(\alpha)}t^{\alpha-1}+\frac{c_2}{\Gamma(\alpha+1)}t^{\alpha} + c_1t^{-1}\sum_{n=1}^{\infty}\frac{(||q||_{\infty}+\Lambda)^n}{\Gamma(2n\alpha+\alpha)}t^{2n\alpha+\alpha}+c_2\sum_{n=1}^{\infty}
\frac{(||q||_{\infty}+\Lambda)^n}{\Gamma(2n\alpha+\alpha+1)}t^{2n\alpha+\alpha}\\
&=c_1t^{-1}\sum_{n=0}^{\infty}\frac{(||q||_{\infty}+\Lambda)^n}{\Gamma(2n\alpha+\alpha)}t^{2n\alpha+\alpha}+c_2\sum_{n=0}^{\infty}
\frac{(||q||_{\infty}+\Lambda)^n}{\Gamma(2n\alpha+\alpha+1)}t^{2n\alpha+\alpha}\\
&=c_1 t^{\alpha-1}E_{2\alpha,\alpha}((||q||_{\infty}+\Lambda)t^{2\alpha})+c_2 t^{\alpha}E_{2\alpha,\alpha+1}((||q||_{\infty}+\Lambda)t^{2\alpha}).
\end{split}
\end{equation}
Note that for a solution $y(t,\lambda)$ of (\ref{IE0}) to be $C([0,1])$, it is necessary and sufficient that $c_1=0$, i.e., $\mathcal{I}^{1-\alpha}_{0^+}y(t)|_{t=0}=0$.
This then proves the global existence of a solution of (\ref{IE0}) on $[\delta,1]$, $\delta>0$, since $q\in L^{\infty}[0,1]$ for given $c_1$ and $c_2$, as defined in \eqref{fcon0}.



From the proof comes the following {\it a-priori} estimate when $c_1=0$, that is, 
\begin{gather}
|y(t, \lambda ) | \leq c_2 {t^\alpha}\left ( \frac{1} {\Gamma(\alpha +1)} + \bigg | E_{2\alpha, \alpha+1} ((||q||_\infty +\Lambda)t^{2\alpha})\bigg | \right )\nonumber\\
 \leq c_2 \left ( \frac{1} {\Gamma(\alpha +1)} + \bigg | E_{2\alpha, \alpha+1} ((||q||_\infty +\Lambda)t^{2\alpha})\bigg | \right )  \nonumber 
\end{gather}
valid for each $t\in [0,1]$ and all $|\lambda|  < \Lambda$.

The previous bound can be made into an absolute constant by taking the sup over all t and  $|\lambda| < \Lambda$. Of course, the bound goes to infinity as $|\lambda| \to \infty$ over non-real values, as it must. Thus, 
\begin{gather}
 |y(t, \lambda ) | \leq c_2 \left (  \frac{1} {\Gamma(\alpha +1)} + \sup_{|\lambda|<\Lambda, t\in [0,1]}  \bigg | E_{2\alpha, \alpha+1} ((||q||_\infty +\Lambda)t^{2\alpha}) \bigg | \right )  \nonumber \\
 = c_2 \left (  \frac{1} {\Gamma(\alpha +1)} +   | E_{2\alpha, \alpha+1} ((||q||_\infty +\Lambda))  | \right ) := c_3.\label{sup}
 \end{gather}
 for all $|\lambda| < \Lambda, t \in [0,1]$. Uniqueness follows easily by means of Gronwall's inequality, as usual. Let $\varepsilon >0$. Assume that  \eqref{IE0} has two solutions $y(t,\lambda), z(t,\lambda)$. Since $q \in L^\infty[0,1]$ and $|\lambda| < \Lambda$ we can derive that,
 $$|y(t,\lambda) - z(t,\lambda)| \leq \varepsilon e^{\frac{1}{\Gamma(2\alpha)} (||q||_\infty + \Lambda) \frac{t^{2\alpha}}{2\alpha}}.$$
 and since $t \in [0,1]$, we get $$|y(t,\lambda) - z(t,\lambda)| \leq \, O(\varepsilon)$$ where the $O$-term can be made independent of both $t, \lambda$. Letting $\varepsilon \to 0$ yields uniqueness for $t \in [0,1]$ and $|\lambda| < \Lambda$.

\section{Another integral equation}
In the previous section we showed that \eqref{IE0} has a solution that, for each $\lambda\in \mathbf{C}$, exists on $[0,1]$, is unique, and is continuous there if and only if $c_1=0$. On the other hand, if $c_1 \neq 0$ then the solution is merely continuous on all compact subsets of $(0,1]$. In this section we find another expression for the integral equation which is equivalent to both \eqref{IE0} and the problem (\ref{fe30}) with boundary conditions (\ref{fcon0}).


\begin{lemma}\label{ft} For $0<\alpha<1$ and $0<t<1$, we have
\[
-^{c}\mathcal{D}^{\alpha}_{0^+}\mathcal{D}^{\alpha}_{0^+}\left(t^{\alpha-1}E_{2\alpha,\alpha}(-\lambda t^{2\alpha})\right) = \lambda t^{\alpha-1}E_{2\alpha,\alpha}(-\lambda t^{2\alpha})
\]
\end{lemma}
\begin{proof}
Using properties of the Mittag-Leffler function we can write
\begin{equation}\label{lem1e1}
\begin{split}
\mathcal{D}^{\alpha}_{0^+}\left(t^{\alpha-1}E_{2\alpha,\alpha}(-\lambda t^{2\alpha})\right)&=\mathcal{D}^{\alpha}_{0^+}\left(
\sum_{k=0}^{\infty}\frac{(-\lambda)^kt^{2\alpha k+\alpha-1}}{\Gamma(2\alpha k+\alpha)}\right)\\
&=\sum_{k=0}^{\infty}\frac{(-\lambda)^k \mathcal{D}^{\alpha}_{0^+}\left(t^{2\alpha k+\alpha-1}\right)}{\Gamma(2\alpha k+\alpha)}\\
&=\sum_{k=0}^{\infty}\frac{(-\lambda)^kt^{2\alpha k-1}}{\Gamma(2\alpha k)}\\
&=t^{-1}E_{2\alpha,0}(\lambda t^{2\alpha})\\
&=-\lambda t^{2\alpha-1}E_{2\alpha,2\alpha}(-\lambda t^{2\alpha}),
\end{split}
\end{equation}
in which the third and the last equalities come from Property \ref{Dt} and Property \ref{rele}, respectively. Now, taking the left Caputo fractional derivative of both sides of (\ref{lem1e1}) we get
\begin{equation}\nonumber
\begin{split}
-^{c}\mathcal{D}^{\alpha}_{0^+}\mathcal{D}^{\alpha}_{0^+}\left(t^{\alpha-1}E_{2\alpha,\alpha}(-\lambda t^{2\alpha})\right)&=
  ^{c}\mathcal{D}^{\alpha}_{0^+}\left(\lambda t^{2\alpha-1}E_{2\alpha,2\alpha}(-\lambda t^{2\alpha})\right)\\
  &=\lambda ^{c}\mathcal{D}^{\alpha}_{0^+}\left(\sum_{k=0}^{\infty}\frac{(-\lambda)^kt^{2\alpha k +2\alpha-1}}
  {\Gamma(2\alpha k+2\alpha)}\right)\\
  &=\lambda \left(\sum_{k=0}^{\infty}\frac{(-\lambda)^k {^{c}\mathcal{D}^{\alpha}_{0^+}}\left(t^{2\alpha k +2\alpha-1}\right)}
  {\Gamma(2\alpha k+2\alpha)}\right)\\
  &= \lambda t^{\alpha-1}\sum_{k=0}^{\infty}\frac{(-\lambda)^kt^{2\alpha k}}{\Gamma(2\alpha k +\alpha)}\\
  &= \lambda t^{\alpha-1}E_{2\alpha,\alpha}(-\lambda t^{2\alpha})
\end{split}
\end{equation}
as required.
\end{proof}

\begin{lemma}\label{st} For $0<\alpha<1$ and $0<t<1$, we have
\[
-^{c}\mathcal{D}^{\alpha}_{0^+}\mathcal{D}^{\alpha}_{0^+}\left(t^{\alpha}E_{2\alpha,\alpha+1}(-\lambda t^{2\alpha})\right) = \lambda t^{\alpha}E_{2\alpha,\alpha+1}(-\lambda t^{2\alpha})
\]
\end{lemma}
\begin{proof}
Once again, using the properties of the Mittag-Leffler function we can write
\begin{equation}\label{lem2e1}
\begin{split}
\mathcal{D}^{\alpha}_{0^+}\left(t^{\alpha}E_{2\alpha,\alpha+1}(-\lambda t^{2\alpha})\right)&=\mathcal{D}^{\alpha}_{0^+}\left(
\sum_{k=0}^{\infty}\frac{(-\lambda)^kt^{2\alpha k+\alpha}}{\Gamma(2\alpha k+\alpha+1)}\right)\\
&=\sum_{k=0}^{\infty}\frac{(-\lambda)^k \mathcal{D}^{\alpha}_{0^+}\left(t^{2\alpha k+\alpha}\right)}{\Gamma(2\alpha k+\alpha+1)}\\
&=\sum_{k=0}^{\infty}\frac{(-\lambda)^kt^{2\alpha k}}{\Gamma(2\alpha k+1)}\\
&=E_{2\alpha,1}(-\lambda t^{2\alpha}).
\end{split}
\end{equation}
in which the third  equality comes from Property \ref{Dtc}. Now, taking the left Caputo fractional derivative of both sides of (\ref{lem2e1}) we get
\begin{equation}\nonumber
\begin{split}
{^c}\mathcal{D}^{\alpha}_{0^+}\mathcal{D}^{\alpha}_{0^+}\left(t^{\alpha}E_{2\alpha,\alpha+1}(-\lambda t^{2\alpha})\right)&=
  {^c}\mathcal{D}^{\alpha}_{0^+}\left(E_{2\alpha,1}(\lambda t^{2\alpha})\right)\\
  &= {^c}\mathcal{D}^{\alpha}_{0^+}\left(\sum_{k=0}^{\infty}\frac{(-\lambda)^kt^{2\alpha k }}
  {\Gamma(2\alpha k+1)}\right)\\
  &=\sum_{k=0}^{\infty}\frac{(-\lambda)^k {^{c}\mathcal{D}^{\alpha}_{0^+}}\left(t^{2\alpha k }\right)}
  {\Gamma(2\alpha k+1)}\\
  &=\sum_{k=1}^{\infty}\frac{(-\lambda)^kt^{\alpha (2k-1)}}{\Gamma(1+\alpha(2k-1))}\\
  &= \sum_{k=0}^{\infty}\frac{(-\lambda)^{k+1}t^{\alpha (2k+1)}}{\Gamma(1+\alpha(2k+1))}\\
  &= -\lambda t^{\alpha}\sum_{k=0}^{\infty}\frac{(-\lambda)^{k}t^{2\alpha k}}{\Gamma(2k\alpha+\alpha+1)}\\
  &=-\lambda t^{\alpha}E_{2\alpha,\alpha+1}(-\lambda t^{2\alpha})
\end{split}
\end{equation}
as desired.
\end{proof}
\begin{lemma}\label{tt} For $0<\alpha<1$ and $0<t<1$, we have
\[
-^{c}\mathcal{D}^{\alpha}_{0^+}\mathcal{D}^{\alpha}_{0^+}\left(\int_0^t(t-s)^{2\alpha-1}E_{2\alpha,2\alpha}
(-\lambda(t-s)^{2\alpha})q(s)y(s)ds\right) = -q(t)y(t)+\lambda \int_0^t(t-s)^{2\alpha-1}E_{2\alpha,2\alpha}
(-\lambda(t-s)^{2\alpha})q(s)y(s)ds
\]
\end{lemma}
\begin{proof}
Let $c_4 = 1/\Gamma(1-\alpha)$. Observe that,
\begin{equation}\label{ialpha}
\begin{split}
  \mathcal{I}^{1-\alpha}_{0^+}\int_0^t(t-s)^{2\alpha-1}E_{2\alpha,2\alpha}
(-\lambda(t-s)^{2\alpha})q(s)y(s)ds&=c_4 \int_0^t\frac{\int_0^r(r-s)^{2\alpha-1}E_{2\alpha,2\alpha}(-\lambda(r-s)^{2\alpha})q(s)y(s)ds}{(t-r)^{\alpha}}dr\\
&=c_4\int_0^t q(s)y(s)\left(  \int_s^t \frac{(r-s)^{2\alpha-1}}{(t-r)^{\alpha}}E_{2\alpha,2\alpha}(-\lambda(r-s)^{2\alpha})dr     \right)ds\\
&=c_4\int_0^t q(s)y(s)\left( \sum_{k=0}^{\infty}\frac{(-\lambda)^k}{\Gamma(2\alpha k +2\alpha)}
\int_s^t \frac{(r-s)^{2\alpha-1+2\alpha k}}{(t-r)^{\alpha}} dr \right) ds\\
&=\int_0^t q(s)y(s)\left( \sum_{k=0}^{\infty}\frac{(-\lambda)^k(t-s)^{2\alpha k + \alpha}}{\Gamma(2\alpha k +\alpha +1)}
  \right) ds\\
 &= \int_0^t q(s)y(s) (t-s)^{\alpha}E_{2\alpha,\alpha+1}(-\lambda(t-s)^{2\alpha})ds.
\end{split}
\end{equation}
Next, differentiating both sides of \eqref{ialpha} with respect to $t$ and noting that $\mathcal{D}^{\alpha}_{0^+}=D(\mathcal{I}^{1-\alpha}_{0^+})$ we find,
\begin{equation}\label{first_part}
\begin{split}
\mathcal{D}^{\alpha}_{0^+}\left(\int_0^t(t-s)^{2\alpha-1}E_{2\alpha,2\alpha}
(-\lambda(t-s)^{2\alpha})q(s)y(s)ds\right)&=\int_0^t(t-s)^{\alpha-1}E_{2\alpha,\alpha}
(-\lambda(t-s)^{2\alpha})q(s)y(s)ds.
\end{split}
\end{equation}
as $\frac{d}{dt}(t^\alpha E_{2\alpha,\alpha+1}(-\lambda t^{2\alpha}))=t^{\alpha-1} E_{2\alpha,\alpha}(-\lambda t^{2\alpha})$.
Next, we are going to take the left Caputo fractional derivative of both sides of (\ref{first_part}). However, since the right hand side of (\ref{first_part}) as a function of $t$ is zero at $t=0$, we can use Property \ref{ricap} and replace the Caputo fractional derivative $^{c}\mathcal{D}^{\alpha}_{0^+}$ by the Riemann-Liouville one $\mathcal{D}^{\alpha}_{0^+}$. In order to do so, first we need to apply $\mathcal{I}^{1-\alpha}_{0^+}$ followed by the classical derivative of the right-hand-side of (\ref{first_part}) as follows,
\begin{equation}\nonumber
\begin{split}
 \mathcal{I}^{1-\alpha}_{0^+}  \int_0^t(t-s)^{\alpha-1}E_{2\alpha,\alpha}
(-\lambda(t-s)^{2\alpha})q(s)y(s)ds &= c_4\int_0^t\frac{\int_0^r(r-s)^{\alpha-1}E_{2\alpha,\alpha}(-\lambda(r-s)^{2\alpha})q(s)y(s)ds}{(t-r)^{\alpha}}dr\\
&=c_4 \int_0^t q(s)y(s)\left(  \int_s^t \frac{(r-s)^{\alpha-1}}{(t-r)^{\alpha}}E_{2\alpha,\alpha}(-\lambda(r-s)^{2\alpha})dr     \right)ds\\
&=c_4 \int_0^t q(s)y(s)\left( \sum_{k=0}^{\infty}\frac{(-\lambda)^k}{\Gamma(2\alpha k +\alpha)}
\int_s^t \frac{(r-s)^{\alpha-1+2\alpha k}}{(t-r)^{\alpha}} dr \right) ds\\
&=\int_0^t q(s)y(s)\left( \sum_{k=0}^{\infty}\frac{(-\lambda)^k(t-s)^{2\alpha k }}{\Gamma(2\alpha k  +1)}
  \right) ds\\
 &= \int_0^t q(s)y(s) E_{2\alpha,1}(-\lambda(t-s)^{2\alpha})ds.
\end{split}
\end{equation}
Taking the derivative of the previous equation and using the fact stated in the previous paragraph, we get
\begin{equation}\label{der}
\begin{split}
 ^{c}\mathcal{D}^{\alpha}_{0^+}  \int_0^t(t-s)^{\alpha-1}E_{2\alpha,\alpha}
(-\lambda(t-s)^{2\alpha})q(s)y(s)ds &=q(t) y(t) + \int_0^tq(s)y(s)(t-s)^{-1}E_{2\alpha,0}(-\lambda(t-s)^{2\alpha})ds\\
&=q(t)y(t) -\lambda \int_0^t(t-s)^{2\alpha-1}E_{2\alpha,2\alpha}(-\lambda(t-s)^{2\alpha})q(s)y(s)ds.
\end{split}
\end{equation}
where we used Property \ref{rele} to arrive at the second equality above. Combining (\ref{first_part}) and (\ref{der}) completes the proof.
\end{proof}
\begin{theorem}\label{TIE}
For $1/2<\alpha<1$, the integral equation 
\begin{gather}\label{IE}
y(t,\lambda) =  {c_1}t^{\alpha-1}\,E_{2\alpha,\alpha}(-\lambda t^{2\alpha})+ c_2 t^{\alpha}E_{2\alpha,\alpha+1}(-\lambda t^{2\alpha})+ \int_0^t(t-s)^{2\alpha-1} E_{2\alpha,2\alpha}(-\lambda (t-s)^{2\alpha}) q(s)\, y(s,\lambda)\, ds 
\end{gather}
satisfies \eqref{fe30} with initial conditions $\mathcal{I}_{0^+}^{1-\alpha}y(t)|_{t=0}=c_1$ and $\mathcal{D}_{0^+}^{\alpha}y(t)|_{t=0}=c_2$ in which $c_1$ and $c_2$ are given constants, and that this solution is unique.
\end{theorem}
\begin{proof}
We apply $-^{c}\mathcal{D}^{\alpha}_{0^+}\mathcal{D}^{\alpha}_{0^+}$ on both sides of (\ref{IE}) to find,
\begin{equation}\label{apply}
\begin{split}
-^{c}\mathcal{D}^{\alpha}_{0^+}\mathcal{D}^{\alpha}_{0^+}\left(    y(t,\lambda)     \right) &=-^{c}\mathcal{D}^{\alpha}_{0^+}\mathcal{D}^{\alpha}_{0^+}\left( {c_1}t^{\alpha-1}\,E_{2\alpha,\alpha}(-\lambda t^{2\alpha})+ c_2 t^{\alpha}E_{2\alpha,\alpha+1}(-\lambda t^{2\alpha})\right ) + \\ 
& - {^c}\mathcal{D}^{\alpha}_{0^+}\mathcal{D}^{\alpha}_{0^+}\left ( \int_0^t(t-s)^{2\alpha-1} E_{2\alpha,2\alpha}(-\lambda (t-s)^{2\alpha}) q(s)\, y(s,\lambda)\, ds  \right) \\
&=\lambda  {c_1}t^{\alpha-1}\,E_{2\alpha,\alpha}(-\lambda t^{2\alpha}) + \lambda c_2 t^{\alpha}E_{2\alpha,\alpha+1}(-\lambda t^{2\alpha})-q(t)y(t)+ \\
&  \lambda \int_0^t(t-s)^{2\alpha-1} E_{2\alpha,2\alpha}(-\lambda (t-s)^{2\alpha}) q(s)\, y(s,\lambda)\, ds\\
&=-q(t)y(t) + \lambda \left(  {c_1}t^{\alpha-1}\,E_{2\alpha,\alpha}(-\lambda t^{2\alpha}) +   c_2 t^{\alpha}E_{2\alpha,\alpha+1}(-\lambda t^{2\alpha})\right )  + \\
& \lambda\, \left (\int_0^t(t-s)^{2\alpha-1} E_{2\alpha,2\alpha}(-\lambda (t-s)^{2\alpha}) q(s)\, y(s,\lambda)\, ds  \right)\\
&=-q(t)y(t) + \lambda \left( y(t,\lambda)\right),
\end{split}
\end{equation}
in which second equality come from Lemma \ref{ft}, Lemma \ref{st}, and Lemma \ref{tt}. We verify the initial conditions. Taking $\mathcal{I}^{1-\alpha}_{0^+}$ of both sides (\ref{IE}), we get,
\begin{equation}\label{first_init}
\begin{split}
\mathcal{I}^{1-\alpha}_{0^+}\left(    y(t,\lambda)     \right) &=\mathcal{I}^{1-\alpha}_{0^+}\left( {c_1}t^{\alpha-1}\,E_{2\alpha,\alpha}(-\lambda t^{2\alpha}) + c_2 t^{\alpha}E_{2\alpha,\alpha+1}(-\lambda t^{2\alpha})\right ) + \\
& \mathcal{I}^{1-\alpha}_{0^+}\left(  \int_0^t(t-s)^{2\alpha-1} E_{2\alpha,2\alpha}(-\lambda (t-s)^{2\alpha}) q(s)\, y(s,\lambda)\, ds  \right) \\
&=  c_1E_{2\alpha,1}(-\lambda t^{2\alpha}) + c_2 tE_{2\alpha,2}(-\lambda t^{2\alpha})+ \int_0^t(t-s)^{\alpha} E_{2\alpha,2\alpha+1}(-\lambda (t-s)^{2\alpha}) q(s)\, y(s,\lambda)\, ds,
\end{split}
\end{equation}
where the third term of the second equality comes from (\ref{ialpha}). Since $E_{2\alpha,1}(-\lambda t^{2\alpha})|_{t=0}=1$ and the other two terms of the above equality vanish when $t=0$, we have verified the first initial condition. Again Taking $\mathcal{D}^{\alpha}_{0^+}$ on both sides (\ref{IE}), we can find,
\begin{equation}\label{second_init}
\begin{split}
\mathcal{D}^{\alpha}_{0^+}\left(    y(t,\lambda)     \right) &=\mathcal{D}^{\alpha}_{0^+}\left( {c_1}t^{\alpha-1}\,E_{2\alpha,\alpha}(-\lambda t^{2\alpha})+ c_2 t^{\alpha}E_{2\alpha,\alpha+1}(-\lambda t^{2\alpha})\right ) + \\
& \mathcal{D}^{\alpha}_{0^+}\left( \int_0^t(t-s)^{2\alpha-1} E_{2\alpha,2\alpha}(-\lambda (t-s)^{2\alpha}) q(s)\, y(s,\lambda)\, ds  \right) \\
&=  -c_1\lambda t^{2\alpha-1} E_{2\alpha,2\alpha}(-\lambda t^{2\alpha}) + c_2 E_{2\alpha,1}(-\lambda t^{2\alpha})+ \int_0^t(t-s)^{\alpha-1} E_{2\alpha,\alpha}(-\lambda (t-s)^{2\alpha}) q(s)\, y(s,\lambda)\, ds,
\end{split}
\end{equation}
where the second equality above comes from (\ref{lem1e1}), (\ref{lem2e1}), and (\ref{first_part}). The second initial condition can readily be obtained by substituting $t=0$ in (\ref{second_init}).
\end{proof}

\section{Analyticity of solutions with respect to the parameter $\lambda$}
In this section we show that the solutions \eqref{IE0} or \eqref{IE} are, generally speaking, entire functions of the parameter $\lambda$ for each $t$ under consideration and $\lambda \in \mathbf{C}.$ First ,we show continuity with respect to said parameter. Consider the case where $c_1=0$, i.e., $y \in C[0,1]$.

\begin{lemma}\label{lemma2}
Let $y \in C[0,1]$, $\lambda \in \mathbf{C}$. Then, for each fixed $t\in [0,1]$, $y(t,\lambda)$ is continuous with respect to $\lambda$ . 
\end{lemma}

\begin{proof}  Let $\Lambda > 0$ be arbitrary but fixed, and let $|\lambda|, |\lambda_0| < \Lambda$. Using \eqref{IE},\begin{equation}\nonumber
\begin{split}
y(t,\lambda)-y(t,\lambda_0)&=\frac{1}{\Gamma(2\alpha)}\int_0^t (t-s)^{2\alpha-1}\left((q(s)-\lambda)y(s,\lambda)-(q(s)-\lambda_0)y(s,\lambda_0)\right)ds\\
&=\frac{1}{\Gamma(2\alpha)}\int_0^t (t-s)^{2\alpha-1}\left((\lambda_0-\lambda)y(s,\lambda)+(q(s)-\lambda_0)(y(s,\lambda)-y(s,\lambda_0))\right)ds.
\end{split}
\end{equation}
So,
\begin{equation}\label{intes}
\begin{split}
y(t,\lambda)-y(t,\lambda_0)=-(\lambda-\lambda_0)\frac{1}{\Gamma(2\alpha)}\int_0^t (t-s)^{2\alpha-1}y(s,\lambda)ds+\frac{1}{\Gamma(2\alpha)}\int_0^t (t-s)^{2\alpha-1}(q(s)-\lambda_0)(y(s,\lambda)-y(s,\lambda_0))ds.
\end{split}
\end{equation}
Now, let $\epsilon>0$ and $|\lambda-\lambda_0|<\delta$ where $\delta>0$ is to be chosen later. Then,
\begin{equation}\nonumber
|y(t,\lambda)-y(t,\lambda_0)|\leq \delta \frac{1}{\Gamma(2\alpha)}\int_0^t(t-s)^{2\alpha-1}|y(s,\lambda)|ds+\frac{1}{\Gamma(2\alpha)}\int_0^t(t-s)^{2\alpha-1}|q(s)-\lambda_0||y(s,\lambda)-y(s,\lambda_0)|ds.
\end{equation}
Using \eqref{sup} and Gronwall's inequality, we get
\begin{equation}\nonumber
\begin{split}
|y(t,\lambda)-y(t,\lambda_0)|&\leq \frac{\delta c_3 t^{2\alpha}}{2\alpha\,\Gamma(2\alpha)} + \frac{1}{\Gamma(2\alpha)}\int_0^t(t-s)^{2\alpha-1}|q(s)-\lambda_0||y(s,\lambda)-y(s,\lambda_0)|ds\\
&\leq \frac{\delta c_3}{\Gamma(2\alpha+1)} \, e^{\frac{1}{\Gamma(2\alpha)}\int_0^1(t-s)^{2\alpha-1}|q(s)-\lambda_0|\,ds}\\
& \leq\, \frac{ \delta c_3}{\Gamma(2\alpha+1)} \, \, e^{\frac{1}{\Gamma(2\alpha)}\int_0^ 1 (1-s)^{2\alpha-1}|q(s)-\lambda_0|\,ds} := C\delta\\
\end{split}
\end{equation}
where 
$$ C=\frac{c_3}{\Gamma(2\alpha+1)} \, \, e^{\frac{1}{\Gamma(2\alpha)}\int_0^ 1 (1-s)^{2\alpha-1}|q(s)-\lambda_0|\,ds}$$
is a function of $\alpha$ and $\lambda_0$ only as $q \in L^\infty(0,1)$. Thus, for any $t\in [0,1]$, the continuity of $y(t,\lambda)$ follows by choosing $\delta<\frac{\varepsilon}{C}$. It also follows from this that,
\begin{equation}\label{bound3}
\sup_{t\in [0,1]} |y(t,\lambda)-y(t,\lambda_0)| < \varepsilon,\quad |\lambda-\lambda_0| < \delta.
\end{equation}
\end{proof}

Next, we consider the differentiability of $y(t,\lambda)$ with respect to $\lambda$. 

\begin{lemma}\label{lemma2x}
Let $y \in C[0,1]$, $\lambda \in \mathbf{C}$. Then, for each fixed $t\in [0,1]$, $y(t,\lambda)$ is differentiable with respect to $\lambda$. 
\end{lemma}

\begin{proof}
As before let $|\lambda| < \Lambda$, $t\in [0,1]$. Equation (\ref{intes}) can be rewritten as 
\begin{equation}\nonumber
\begin{split}
\frac{y(t,\lambda)-y(t,\lambda_0)}{\lambda-\lambda_0}=-\frac{1}{\Gamma(2\alpha)}\int_0^t (t-s)^{2\alpha-1}y(s,\lambda)\, ds+\frac{1}{\Gamma(2\alpha)}\int_0^t (t-s)^{2\alpha-1}(q(s)-\lambda_0)\frac{y(s,\lambda)-y(s,\lambda_0)}{\lambda-\lambda_0}\,ds.
\end{split}
\end{equation} 
As $y(t,\lambda_0)$ is given, we define $h(t,\lambda_0)$ to be the unique solution of the Volterra integral equation of the second kind,
\begin{equation}\nonumber
\begin{split}
h(t,\lambda_0)=-\frac{1}{\Gamma(2\alpha)}\int_0^t (t-s)^{2\alpha-1}y(s,\lambda_0)\,ds+\frac{1}{\Gamma(2\alpha)}\int_0^t (t-s)^{2\alpha-1}(q(s)-\lambda_0)h(s,\lambda_0)\,ds.
\end{split}
\end{equation} 
So,
\begin{equation}\nonumber
\begin{split}
\bigg |\frac{y(t,\lambda)-y(t,\lambda_0)}{\lambda-\lambda_0}-h(t,\lambda_0)\bigg |&\leq \frac{1}{\Gamma(2\alpha)}\int_0^t (t-s)^{2\alpha-1}|y(s,\lambda)-y(s,\lambda_0) |ds\\
&+\frac{1}{\Gamma(2\alpha)}\int_0^t (t-s)^{2\alpha-1}|q(s)-\lambda_0|\bigg |\frac{y(s,\lambda)-y(s,\lambda_0)}{\lambda-\lambda_0}-h(s,\lambda_0)\bigg |ds.
\end{split}
\end{equation} 
Let $\varepsilon>0$ and choose $\delta>0$ as in \eqref{bound3}. Using Gronwall's inequality and \eqref{bound3} we get, for $t \in [0,1]$,
\begin{equation}\label{diff}
\begin{split}
\bigg |\frac{y(t,\lambda)-y(t,\lambda_0)}{\lambda-\lambda_0}-h(t,\lambda_0)\bigg |&\leq \frac{\varepsilon}{2\alpha \Gamma(2\alpha)} +\frac{1}{\Gamma(2\alpha)}\int_0^t (t-s)^{2\alpha-1}|q(s)-\lambda_0|\bigg |\frac{y(s,\lambda)-y(s,\lambda_0)}{\lambda-\lambda_0}-h(s,\lambda_0)\bigg |ds,\\
&\leq \frac{\varepsilon}{\Gamma(2\alpha+1)} \, e^{\frac{1}{\Gamma(2\alpha)}\int_0^1 (t-s)^{2\alpha-1}|q(s)-\lambda_0|\,ds} = O(\varepsilon).
\end{split}
\end{equation}
for $\lambda$ near $\lambda_0$ since, for $t \in [0,1]$, $\int_0^1 (t-s)^{2\alpha-1}|q(s)-\lambda_0|\,ds=O(1)$. Thus, 
\[
\frac{\partial y(t,\lambda)}{\partial \lambda}|_{\lambda=\lambda_0}:=\lim_{\lambda\rightarrow\lambda_0}\frac{y(t,\lambda)-y(t,\lambda_0)}{\lambda-\lambda_0}=h(t,\lambda_0),
\]
exists at $\lambda_0$. Since $\lambda_0$ is arbitrary $y_\lambda(t,\lambda)$ exists for all $\lambda$ with $|\lambda| < \Lambda$, real or complex and the result follows. 
\end{proof}

\begin{theorem}\label{entire}
For each $t\in [0,1]$, $y(t,\lambda)$ is an entire function of $\lambda$.
\end{theorem}

\begin{proof} This follows from Lemma~\ref{lemma2x} since $\lambda \in \mathbf{C}$ and $|\lambda | < \Lambda$ where $\Lambda>0$ is arbitrary.
\end{proof}

\section{A Dirichlet type problem}

Let $y \in C[0,1]$, $\lambda \in \mathbf{C}$ be fixed. In this case we note that the first of the boundary conditions \eqref{fcon0} is equivalent to the usual fixed end (Dirichlet) boundary conditions, that is,
$$ y \in C[0,1] \quad \iff \mathcal{I}_{0^+}^{1-\alpha}y(t,\lambda)|_{t=0}=0\quad  \iff \quad y(0,\lambda)=0.$$
For the continuity assumption implies that there is a number $M$ such that $|y(t,\lambda)|\leq M$, for all $t\in [0,1]$. Thus,
$$|\mathcal{I}_{0^+}^{1-\alpha}y(t,\lambda) | \leq \frac{M}{\Gamma(\alpha)} \int_0^t (t-s)^{-\alpha}\, ds = \frac{M\,t^{1-\alpha}}{(1-\alpha)\Gamma(\alpha)}, $$ and so $\mathcal{I}_{0^+}^{1-\alpha}y(t,\lambda)|_{t=0}=0$. On the other hand \eqref{IE0} now implies that $c_1=0$, i.e., $y(0,\lambda)=0$, so that $y \in C[0,1].$ However, the condition $y(1,\lambda)=0$ is independent of the statement that $\mathcal{I}_{0^+}^{1-\alpha}y(t,\lambda)|_{t=1}=0$.

Since, for any $z\neq 0$,  the Mittag-Leffler functions satisfy $$E_{\delta,\delta}(z)=\frac{1}{z} E_{\delta,0}(z) ,$$ we get
\begin{equation}\label{equl}
t^{2\alpha-1}E_{2\alpha,2\alpha}(-\lambda t^{2\alpha})= - \frac{1}{\lambda t}\,E_{2\alpha,0}(-\lambda t^{2\alpha}) .
\end{equation}
Hence, using  (\ref{IE}) and (\ref{equl}) we get 
\begin{equation}\label{volt}
y(t,\lambda)={c_1}\,{t^{\alpha - 1}}\,E_{2\alpha,\alpha}(-\lambda t^{2\alpha})+ c_2\, t^{\alpha}E_{2\alpha,\alpha+1}(-\lambda t^{2\alpha}) - \int_0^t \frac{E_{2\alpha,0}(-\lambda (t-s)^{2\alpha})}{\lambda (t-s)}q(s)y(s,\lambda)\,ds.
\end{equation}

{\bf Remark 1:}\quad  When $\alpha\to 1$, the integral equation (\ref{volt}) becomes  
\begin{equation}\label{cvolt}
y(t,\lambda)=y(0,\lambda)\cos(\sqrt{\lambda}t)+y^\prime(0,\lambda)\frac{\sin(\sqrt{\lambda}t)}{\sqrt{\lambda}}+\int_0^t \frac{\sin(\sqrt{\lambda}(t-s))}{\sqrt{\lambda}}q(s)y(s,\lambda)ds,
\end{equation}
which is exactly the integral equation equivalent of  the classical Sturm-Liouville equation $-y^{\prime\prime}+q(t)y=\lambda y$ for $\lambda >0$.

{\bf Remark 2:}\quad Observe that, for each $\alpha$, 
\[ \lim_{s\to t^-} - \frac{E_{2\alpha,0}(-\lambda (t-s)^{2\alpha})}{\lambda (t-s)} = \displaystyle \left \{ \begin{array}{ll}
0, & \mbox{if \ \ $\alpha \in (1/2 , 1]$},\\
\displaystyle  1, & \mbox{if \ \ $\alpha =1/2$}.\\
\end{array}
\right.
\]

 and so, for each $1/2 < \alpha < 1$, the kernel appearing in \eqref{volt} is uniformly bounded on $[0,1]$. This agrees with the equivalent result for the classical case (\ref{cvolt}).
 
 

\section{Existence and asymptotic distribution of the eigenvalues}
Without loss of generality we may assume that $c_2=1$ in \eqref{volt} and $y(t,\lambda)$ is the corresponding solution. In the sequel we always assume that $1/2 < \alpha < 1$.

\begin{lemma}\label{lemma02}
For each $t \in [0,1]$, $1/2<\alpha<1$, and $|\arg(-\lambda)|\leq \mu$ where $\mu \in (\alpha\pi,\pi)$, we have\\ $\left |t^{\alpha}E_{2\alpha,\alpha+1}(-\lambda t^{2\alpha})\right | \to 0$ as $|\lambda|\to \infty$.
\end{lemma}
\begin{proof}
By (\ref{asyx1}) we can write
\begin{equation}\nonumber
\begin{split}
t^{\alpha}E_{2\alpha,\alpha+1}(-\lambda t^{2\alpha}) & = t^{\alpha}\left( \frac{1}{2\alpha}(-\lambda t^{2\alpha})^{\frac{1-(\alpha+1)}{2\alpha}}\right)\exp \bigg \{  (-\lambda t^{2\alpha})^{\frac{1}{2\alpha}}\bigg \} + \text{\Large {O}}\left(  \frac{1}{\lambda}\right) \\
& =  -\frac{i}{2\alpha\sqrt{\lambda}}\exp \bigg \{  (-\lambda)^{\frac{1}{2\alpha}}t \bigg \} + \text{\Large {O}}\left(  \frac{1}{\lambda}\right)\\
& =  -\frac{i}{2\alpha\sqrt{\lambda}}\exp \bigg \{  |\lambda|^{\frac{1}{2\alpha}}\left(\cos(\frac{\arg(-\lambda)}{2\alpha})+i\sin(\frac{\arg(-\lambda)}{2\alpha})\right) t\bigg \} + \text{\Large {O}}\left(  \frac{1}{\lambda}\right).
\end{split}
\end{equation}
Therefore,
\begin{equation}\nonumber
\begin{split}
\left |t^{\alpha}E_{2\alpha,\alpha+1}(-\lambda t^{2\alpha}) \right  |& = \frac{1}{2\alpha\sqrt{\lambda}}\exp \bigg\{  |\lambda|^{\frac{1}{2\alpha}}\cos(\frac{\arg(-\lambda)}{2\alpha})t\bigg\}.
\end{split}
\end{equation}
Regarding the assumption on $\arg(-\lambda)$, we have $\cos(\frac{\arg(-\lambda)}{2\alpha})<0$ and it completes the proof.
\end{proof}


\begin{lemma}\label{lemma03}
For each $t \in [0,1]$, $s \in [0,t]$, $1/2<\alpha<1$, and $|\arg(-\lambda)|\leq \mu$ where $\mu \in (\alpha\pi,\pi)$, we have $\left |\frac{E_{2\alpha,0}(-\lambda (t-s)^{2\alpha})}{\lambda(t-s)}\right |\to 0$ as $|\lambda|\to \infty$.
\end{lemma}

\begin{proof}
By (\ref{asyx1}) we can write
\begin{equation}\nonumber
\begin{split}
\frac{E_{2\alpha,0}(-\lambda (t-s)^{2\alpha})}{\lambda(t-s)} & =\frac{ \left( \frac{1}{2\alpha}(-\lambda (t-s)^{2\alpha})^{\frac{1}{2\alpha}}\right)\exp \bigg \{  (-\lambda (t-s)^{2\alpha})^{\frac{1}{2\alpha}}\bigg \} + \text{\Large {O}}\left(  \frac{1}{\lambda}\right)}{\lambda(t-s)} \\
& = \frac{1}{2\alpha}\frac{(-\lambda)^{\frac{1}{2\alpha}}}{\lambda}\exp \bigg \{  (-\lambda)^{\frac{1}{2\alpha}}(t-s) \bigg \} + \text{\Large {O}}\left(  \frac{1}{\lambda^2}\right).
\end{split}
\end{equation}
Then,

$$\left | \frac{E_{2\alpha,0}(-\lambda (t-s)^{2\alpha})}{\lambda (t-s)}\right |  = \frac{1}{2\alpha |\lambda|^{(2\alpha-1)/2\alpha}} \,  \exp \bigg \{ ( t-s)  |\lambda|^{1/2\alpha} \cos \left (  \frac{\arg(-\lambda)} {2\alpha}  \right )\bigg  \} + \Large {O}\left ( \frac{1}{|\lambda|^2} \right). $$ 
Arguing as in the previous lemma we reach the desired conclusion.

\end{proof}
%

\begin{lemma}\label{lemma05}
For each $t \in [0,1]$, $s \in [0,t]$, $1/2<\alpha<1$, and $|\arg(-\lambda)|\leq \mu$ where $\mu \in (\alpha\pi,\pi)$, we have $\left |\frac{E_{2\alpha,1-\alpha}(-\lambda (t-s)^{2\alpha})}{\lambda (t-s)^{\alpha}}\right |\to 0 $ as $|\lambda|\to \infty$


\end{lemma}
{\begin{proof}
By (\ref{asyx1}) we can write
\begin{equation}\nonumber
\begin{split}
\frac{E_{2\alpha,1-\alpha}(-\lambda (t-s)^{2\alpha})}{\lambda(t-s)^{\alpha}} & =\frac{ \left( \frac{1}{2\alpha}(-\lambda (t-s)^{2\alpha})^{\frac{1-(1-\alpha)}{2\alpha}}\right)\exp \bigg \{  (-\lambda (t-s)^{2\alpha})^{\frac{1}{2\alpha}}\bigg \} + \text{\Large {O}}\left(  \frac{1}{\lambda}\right)}{\lambda(t-s)^{\alpha}} \\
& = \frac{i}{\sqrt{\lambda}2\alpha}\exp \bigg \{  (-\lambda)^{\frac{1}{2\alpha}}(t-s) \bigg \} + \text{\Large {O}}\left(  \frac{1}{\lambda^2}\right).
\end{split}
\end{equation}
Then,
\begin{gather}\label{eq73}
\left | \frac{E_{2\alpha,1-\alpha}(-\lambda (t-s)^{2\alpha})}{\lambda (t-s)^{\alpha}} \right | = \frac{1}{2\alpha \sqrt{|\lambda|}} \,  \exp \bigg \{ ( t-s)  |\lambda|^{1/2\alpha} \cos \left (  \frac{\arg(-\lambda)} {2\alpha}  \right )\bigg  \} .
\end{gather}
The result follows since the exponential term is uniformly bounded. 

\end{proof}

\begin{lemma}\label{lemma04}
For each $t \in [0,1]$, and $1/2< \alpha < 1$, the solution $y(t,\lambda)$ is an entire function of $\lambda$ of order at most $1/2\alpha$.
\end{lemma}
\begin{proof}  Let $\lambda \in \mathbf{C}$. Define $f$ by 
\begin{gather}\label{eqf} y(t,\lambda) = \exp \bigg \{ t |\lambda|^{1/2\alpha} \cos \left (  \frac{\arg(-\lambda)} {2\alpha}  \right )\bigg  \} f(t).
\end{gather}
 Then, using \eqref{volt}, 
\begin{gather*}
f(t) = t^{\alpha}E_{2\alpha,\alpha+1}(-\lambda t^{2\alpha}) \exp \bigg \{ - t |\lambda|^{1/2\alpha} \cos \left (  \frac{\arg(-\lambda)} {2\alpha}  \right )\bigg  \}  \nonumber \\
- \int_0^t \frac{E_{2\alpha,0}(-\lambda (t-s)^{2\alpha})}{\lambda (t-s)} \exp \bigg \{ -( t-s)  |\lambda|^{1/2\alpha} \cos \left (  \frac{\arg(-\lambda)} {2\alpha}  \right )\bigg  \} \ q(s) f(s) \,ds
\end{gather*}

Applying Lemma~\ref{lemma02} there exists $\Lambda \in \mathbf{R^+}$ such that for all $|\lambda| > \Lambda$ we have 
\begin{gather*}
|f(t)| \leq 1+  \frac{1}{2\alpha |\lambda|^{(2\alpha-1)/2\alpha}}\int_0^t |q(s)|\,|f(s)|\, ds\\
\end{gather*}
which, on account of Gronwall's inequality, gives us
\begin{gather}
|f(t)| \leq \exp\left \{ \frac{1}{2\alpha |\lambda|^{(2\alpha-1)/2\alpha}}\int_0^1 |q(s)|\,ds\right \}
\end{gather}
for all sufficiently large $|\lambda|$. Thus, $f \in L^\infty[0,1]$ so that \eqref{eqf} yields, for some $M$,
$$|y(t,\lambda)| \leq M \exp \left ( |\lambda|^{1/2\alpha} \right )$$
and the order claim is verified.
\end{proof}

\begin{lemma}\label{lemma01}
For each $t \in [0,1]$, $\mathcal{I}_{0^+}^{1-\alpha}y(t,\lambda)$ is an entire function of $\lambda$ of order at most $2\alpha$.
\end{lemma}

\begin{proof}
This is clear from the definition, the possible values of $\alpha$, and since $y(t, \lambda)$ is itself entire and of order at most $1/2\alpha$, from Lemma~\ref{lemma04}.
\end{proof}

\begin{lemma}\label{lemma001}
The boundary value problem \eqref{fe30}-\eqref{fcon0} has infinitely many complex eigenvalues (real eigenvalues are not to be excluded here).
\end{lemma}

\begin{proof}
By Lemma~\ref{lemma01}, we know that $\mathcal{I}_{0^+}^{1-\alpha}y(t,\,\lambda)$ is entire for each $t\in [0,1]$, and $1/2 <\alpha<1$  as well. So, the eigenvalues of our problem are given by the zeros of $\mathcal{I}_{0^+}^{1-\alpha}y(1,\,\lambda)$, which must be countably infinite in number since the latter function is of fractional order $1/2\alpha$ (on account of the restriction on $\alpha$). This gives us the existence of infinitely many eigenvalues, generally in $\mathbf{C}$.  
\end{proof}
Next, we give the asymptotic distribution of these eigenvalues when $\alpha$ is either very close to $1/2$ from the right or very close to $1$ from the left. Recall \eqref{volt} with $c_2=1$, so that
\begin{equation}\label{asymp}
y(t,\lambda) =  t^{\alpha}E_{2\alpha,\alpha+1}(-\lambda t^{2\alpha}) -  \int_0^t \frac{ E_{2\alpha,0}(-\lambda (t-s)^{2\alpha})}{\lambda (t-s)}q(s)y(s,\lambda)\,ds.
\end{equation}
Keeping in mind the boundary condition \eqref{fcon0} at $t=1$, we calculate $\mathcal{I}_{0^+}^{1-\alpha}y(t,\,\lambda)$ and then evaluate this at $t=1$ in order to find the dispersion relation for the eigenvalues. 

A straightforward though lengthy calculation using \eqref{asymp} and the definition of the Mittag-Leffler functions show that
\begin{gather}
\mathcal{I}_{0^+}^{1-\alpha}y(t,\,\lambda) = \mathcal{I}_{0^+}^{1-\alpha}\{t^{\alpha}E_{2\alpha,\alpha+1}(-\lambda t^{2\alpha})\} + \mathcal{I}_{0^+}^{1-\alpha} \left (\int_0^t \frac{ E_{2\alpha,0}(-\lambda (t-s)^{2\alpha})}{\lambda (t-s)}q(s)y(s,\lambda)\,ds.
 \right ),\\
 = tE_{2\alpha,2}(-\lambda t^{2\alpha}) +  \frac{1}{\lambda}\int_0^t \frac{ E_{2\alpha,1-\alpha}(-\lambda (t-s)^{2\alpha})}{(t-s)^\alpha }q(s)y(s,\lambda)\,ds
 \end{gather}
 so that the eigenvalues of \eqref{fe30}-\eqref{fcon0} are given by those $\lambda \in \mathbf{C}$ such that
 \begin{gather}\label{disp}
  E_{2\alpha,2}(-\lambda) +  \frac{1}{\lambda}\int_0^1\frac{ E_{2\alpha,1-\alpha}(-\lambda (1-s)^{2\alpha})}{(1-s)^\alpha }q(s)y(s,\lambda)\,ds = 0.
  \end{gather}
Let us consider first the case where $\lambda \in \mathbf{R}$. Lemma~\ref{lemma05} implies that the right side of \eqref{eq73} tends to $0$ as $\lambda \to \infty$. Indeed this, combined with \eqref{eqf}, implies that  
 $$\left | \frac{E_{2\alpha,1-\alpha}(-\lambda (t-s)^{2\alpha})}{\lambda (t-s)^{\alpha}} \, y(s,\lambda) \right | = \Large{O}  \left ( \frac{1}{ \sqrt{|\lambda|}} \right )$$
for all sufficiently large $\lambda$.

 
 Thus, the real eigenvalues of the problem (\ref{fe30})-\eqref{fcon0} become the zeros of a transcendental equation of the form,
\[
E_{2\alpha,2}(-\lambda)+ O\left ( \frac{1}{\sqrt{\lambda}}\right )=0.
\]
We are concerned with the asymptotic behaviour of these real zeros.  Recall  the distribution of the real zeros of $E_{2\alpha,2}(-\lambda)$ in \cite{dm1}. There we showed that, for each $n=0,1,2, \ldots,N^*-1$, where $N^*$ depends on $\alpha$, the interval 
\begin{equation}\label{neg}
I_n(\alpha):=\left(\left(\frac{(2n+\frac{1}{2}+\frac{1}{2\alpha})\pi}{\sin(\frac{\pi}{2\alpha})}\right)^{2\alpha},\quad \left(\frac{(2n+\frac{3}{2}+\frac{1}{2\alpha})\pi}{\sin(\frac{\pi}{2\alpha})}\right)^{2\alpha}\right),
\end{equation}
always contains at least two real zeros of $E_{2\alpha,2}(-\lambda)$. For $\alpha\to 1$, these intervals approach the intervals 
$$ \left ( (2n+1)^2\,\pi^2, (2n+2)^2\,\pi^2 \right ),$$
whose end-points are each eigenvalues of the Dirichlet problem for the classical equation $-y^{\prime\prime} = \lambda\,y$ on $[0,1]$. Since each interval $I_n$ contains two zeros we can denote the first of these two zeros by $\lambda_{2n}(\alpha)$. Equation \eqref{neg} now gives the {\it a-priori} estimate
\begin{equation}\label{neg2}
\left(\frac{(2n+\frac{1}{2}+\frac{1}{2\alpha})\pi}{\sin(\frac{\pi}{2\alpha})}\right)^{2\alpha} \leq   \lambda_{2n}(\alpha)    \leq \left(\frac{(2n+\frac{3}{2}+\frac{1}{2\alpha})\pi}{\sin(\frac{\pi}{2\alpha})}\right)^{2\alpha}.
\end{equation}

For each $\alpha<1$, and close to $1$, and for large $\lambda$, the real zeros of the preceding equation approach those of $E_{2\alpha,2}(-\lambda)$ and spread out towards the end-points of intervals of the form \eqref{neg}. For $\alpha$ close to $1/2$ there are no zeros, the first two zeros appearing only when $\alpha \approx 0.7325$. For $\alpha$ larger than this critical value, the zeros appear in pairs and in intervals of the form \eqref{neg}.  

Next, recall that for  $\alpha <1$ there are  only {\it finitely many} such real zeros, (see \cite{dm1})  their number growing without bound as $\alpha \to 1$. It also follows from Lemma~\ref{lemma001}  that, for each $\alpha$, the remaining infinitely many eigenvalues must be non-real. As $\alpha \to 1^-$ these non-real eigenvalues tend to the real axis thereby forming more and more real eigenvalues until the spectrum is totally real when $\alpha =1$ and the problem then reduces to a (classical) regular Sturm-Liouville problem.

Finally, for $\alpha$ close to $1$, \eqref{neg2}  leads to the approximation,
$$\lambda_{2n}(\alpha)    \approx \left(\frac{(2n+2)\pi}{\sin(\frac{\pi}{2\alpha})}\right)^{2\alpha},$$
from which this, in conjunction with \eqref{neg} and $\alpha \to 1$, we can derive the classical eigenvalue asymptotics, $\lambda_n \sim n^2\pi^2$ as $n \to \infty$.

\section{Closing remarks}

We have shown that the fractional eigenvalue problem 
\begin{equation*}
-^{c}\mathcal{D}_{0^+}^{\alpha}\circ   \mathcal{D}^{\alpha}_{0^+}y(t)+q(t)y(t)=\lambda y(t),\qquad 1/2<\alpha<1,\quad 0 \leq t \leq 1,
\end{equation*}
with mixed Caputo and Riemann-Liouville derivatives subject to the boundary conditions involving the Riemann-Liouville integrals,
\begin{equation*}
\mathcal{I}_{0^+}^{1-\alpha}y(t)|_{t=0}=0,\quad \text{and}\quad \mathcal{I}_{0^+}^{1-\alpha}y(t)|_{t=1}=0,
\end{equation*}
admits, for each $\alpha$ under consideration, and for eigenfunctions that are in $C[0,1]$, a finite number of real eigenvalues and an infinite number of non-real eigenvalues. The real eigenvalues, though finite in number for each $\alpha$, are approximated by \eqref{neg} and \eqref{neg2}, which as $\alpha \to 1$ gives the classical asymptotic relation 
$\lambda_n \sim n^2\pi^2$ as $n \to \infty$.

As $\alpha \to 1^-$ we observe that the spectrum obtained approaches the Sturm-Liouville spectrum of the classical problem 
$$-y^{\prime\prime}+ q(t) y = \lambda y,\quad\quad y(0)=y(1)=0.$$
The same results hold if the eigenfunctions are merely $C(0,1]$ (i.e., $c_1\neq 0$) except that now the latter have an infinite discontinuity at $t=0$ for each $\alpha$. The proofs are identical and are therefore omitted.

\section*{References} 


\end{document}